\documentclass[11pt,a4paper]{amsart}

\usepackage{xy}
\usepackage{amssymb}
\usepackage{color}

\DeclareMathOperator{\End}{End}
\DeclareMathOperator{\Ext}{Ext}
\DeclareMathOperator{\Hom}{Hom}
\DeclareMathOperator{\matring}{M}
\DeclareMathOperator{\rk}{rk}

\newcommand{\ra}{\rightarrow}
\newcommand{\NN}{ \mathbb  N}
\newcommand{\ZZ}{\mathbb Z}
\newcommand{\al}{\alpha}
\newcommand{\La}{\Lambda}
\renewcommand{\mod}{\operatorname{mod}}

\numberwithin{equation}{section}

\xyoption{arrow}
\xyoption{matrix}
\xyoption{frame}

\xyoption{ps}
\xyoption{dvips}

\theoremstyle{plain}
\newtheorem{Thm}{Theorem}
\newtheorem{Prop}{Proposition}[section]

\theoremstyle{definition}

\newtheorem{Number}[Prop]{}

\advance\textheight by 1.5cm
\advance\textwidth by 1.3cm
\advance\evensidemargin by -0.8cm
\advance\oddsidemargin by -0.8cm
\advance\topmargin by -1.0cm

\begin{document}

\title[Extended Dynkin]{Indecomposable representations \\
  for extended Dynkin quivers}

\author{Dirk Kussin}
\author{Hagen Meltzer}

\address{Institut f\"ur Mathematik, Universit\"at Paderborn, 33095
  Pader\-born, Germany}

\email{dirk@math.upb.de}

\address{Instytut Matematyki, Uniwersytet Szczeci\'nski, 70451
  Szczecin, Poland}

\email{meltzer@wmf.univ.szczecin.pl}

\subjclass[2000]{16G20}

\keywords{extended Dynkin quiver, Euclidean quiver, tame quiver,
canonical algebra, domestic type,
  preprojective module, exceptional module, tilting module.}

\thanks{The second author was supported by the Polish Scientific Grant
  KBN 1 P03A 007 27.}

\begin{abstract}
  We describe a method for an explicit determination of indecomposable
  preprojective and preinjective representations for extended Dynkin
  quivers $\Gamma$ over an arbitrary field $K$ by vector spaces and matrices.
  This method uses tilting theory and the explicit knowledge of
  indecomposable modules over the corresponding canonical algebra of
  domestic type. Further, if $K$ is algebraically closed we obtain
  all indecomposable representations for $\Gamma$.
   For the case that $\Gamma$ is
  of type $\widetilde{D}_n$, $ n \geq 4$, with a fixed orientation, we
  determine all indecomposable preprojective representations.
  Moreover, in the case $\widetilde{E}_6$ we present the most complicated
  indecomposable preprojective representations of rank $3$.
\end{abstract}

\maketitle

\section{Introduction}
Let $K$ be a field, $\Gamma$ a quiver and $A=K\Gamma/I$ a finite-dimensional
algebra of quiver type.
 One of the problems in representation theory
is to give normal forms for the indecomposable finite-dimensional left
$A$-modules.  Such a module is given by choosing a finite-dimensional
vector space for each vertex of and a linear map for each arrow of the
quiver such that the relations of the ideal $I$ are satisfied.

The problem to determine all indecomposable modules for an algebra
\emph{explicitly\/} by vector spaces and matrices is in general
difficult. This problem is solved only in very few cases. In
particular, Gabriel computed in 1972 the indecomposable
representations for Dynkin quivers \cite{gabriel:72}.

Concerning \emph{extended\/} Dynkin quivers only partial results are
known.  Already in 1890 Kronecker~\cite{Kronecker} classified pairs of
$n \times m$-matrices up to simultaneous equivalence, solving a
problem raised by Weierstra\ss. In modern terminology this means to
describe the finite dimensional modules over the Kronecker algebra.
For the case of $\widetilde{D}_4$ with subspace orientation, the so
called 4-subspace problem, indecomposable objects were described by
Nazarova~\cite{nazarova:67} and Gelfand and
Ponomarev~\cite{gelfand:ponomarev}.

General information about the structure of the module category of a
path algebra of an extended Dynkin quiver were obtained by Donovan and
Freislich~\cite{donovan:freislich}, and Nazarova~\cite{nazarova:73}.
For the characterization of the regular modules and for historical
remarks we refer to \cite[chapter~11]{gabriel:roiter}.

We recall also that Ringel~\cite{ringel:98} has shown that for every
finite dimensional path algebra $A=K\Gamma$ for a quiver $\Gamma$ each
exceptional module
can be exhibited by matrices containing as coefficients only $0$ and
$1$.  For a path algebra over an extended Dynkin quiver each
indecomposable preprojective (respectively preinjective) module is
exceptional, however explicit descriptions for those modules were not
given in this case.

In this paper we discuss a method for the description of
indecomposable representations of extended Dynkin quivers using our
explicit description of indecomposable modules over domestic canonical
algebras given in \cite{kussin:meltzer:04} and \cite{komoda:meltzer}.
We apply tilting theory which was developed in \cite{brenner:butler}
and \cite{happel:ringel}.  We exploit the fact that for each path
algebra $A=K\Gamma$ of an extended Dynkin quiver $\Gamma$ there is a canonical
algebra of domestic type $\Lambda$ and a tilting module $T$ over
$\Lambda$ such that $\End (T) \simeq A^{op} $ \cite{ringel:84} (see
also \cite{huebner} and \cite[Proposition 6.5]{lenzing:reiten}).
  The
theorem of Brenner and Butler \cite{brenner:butler} ensures that
applying the functor $\Hom_{\La}(T,-)$ to an indecomposable
preprojective left $\Lambda$-module $M$ satisfying $\Ext^1_{\La}
(T,M)=0$ we obtain an indecomposable right $A^{op}$-module, thus an
indecomposable representations $N$ of $\Gamma$.  Moreover, in this way we
 obtain all preprojective indecomposable representations of $\Gamma$. In
this paper we concentrate on the description of the
preprojectives. However we remark that the regular modules can be
treated in the same way, using~\cite[Chapter 4]{kussin:meltzer:04}.
Finally, the indecomposable preinjective representations for
$\Gamma$ can be obtained by duality, i.e.\ by choosing the opposite
orientation of the quiver.


\section{Tilting from domestic canonical algebras to path algebras
of extended Dynkin quivers }\label{sec:tilting}

Canonical algebras were introduced by Ringel in 1984 \cite{ringel:84}
and play an important role in representation theory.
A domestic canonical algebra of quiver type $\Lambda$ is isomorphic to
the path algebra of the quiver

$$
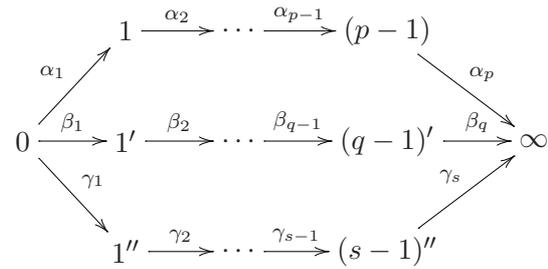
\xy\xymatrixcolsep{2pc}\xymatrix{ & 1 \ar @{->}[r]^-{\alpha_2} &
  \cdots \ar @{->}[r]^-{\alpha_{ p-1 }} &( p-1) \ar
  @{->}[dr]^-{\alpha_{p}} & \\
  0 \ar @{->}[ur]^-{\alpha_1} \ar @{->}[dr]^-{\gamma_1} \ar
  @{->}[r]^-{\beta_1} &
  1' \ar @{->}[r]^-{\beta_2} & \cdots \ar @{->}[r]^-{\beta_{q-1}} &
 ( q-1)' \ar @{->}[r]^-{\beta_q} & \infty \\ & 1''
 \ar @{->}[r]^-{\gamma_2} &
  \cdots \ar @{->}[r]^-{\gamma_{s-1}} &
 (s-1)'' \ar @{->}[ur]^-{\gamma_s} &} \endxy
$$

\noindent
modulo the relation $\gamma_s \dots \gamma_1=\alpha_p \dots \alpha_1 +
\beta_q \dots \beta_1$, where $p$, $q$, $s$ is the length of the upper
(middle, lower, respectively) arm, and where moreover the triple
$(p,q,s)$ is given by $(p,q,1)$ (where $p$, $q\geq 1$), $(p,2,2)$
(where $p\geq 2$), $(3,3,2)$, $(4,3,2)$ or $(5,3,2)$.

Therefore a finite-dimensional left $\Lambda$-module $M$ consists of
finite-dimensional vector spaces $M (i)$ for each point $i$ of the
quiver, and a linear map $M(\alpha)$ for each arrow $\alpha=\alpha_i$,
$\beta_j$ and $\gamma_k$, satisfying the relation
$$M(\gamma_s)\circ\ldots\circ M(\gamma_1)=M(\alpha_p)\circ\ldots\circ
M(\alpha_1) +M(\beta_q)\circ\ldots\circ M(\beta_1).$$
The number $\rk (M)=\dim M(\infty)-\dim M (0)$
is called the
\emph{rank\/} of $M$.
 Then an indecomposable module of positive rank
(negative rank, rank zero, respectively) is preprojective
(preinjective, regular, respectively).
The global structure of the module category looks as follows: There is
precisely one preprojective component and precisely one preinjective
component and the indecomposable regular modules form
tubes~\cite{ringel:84}.

Indecomposable preprojective left modules over $\Lambda$ were
described by explicit matrices in \cite{kussin:meltzer:04} in case
that the characteristic of $K$ is different from $2$ and in
\cite{komoda:meltzer} for an arbitrary field, the last is relevant
only for modules of rank $6$ in the domestic situation $(5,3,2)$.  The
indecomposable $\Lambda$-modules appear in series and are constructed
using a general principle by applying the so called method of
enlargement of matrices and adding identities.  We recall this general
principle and provide for this  the following notations.

Let $n$ and $i$ be natural numbers. Let $I_n$ be the $n\times
n$-identity matrix. Define $$X_{n+i}^{n}=
\begin{bmatrix}
 &  I_n & \\
 \hline
 0 & \cdots & 0\\
 \vdots & & \vdots  \\
 0 & \cdots & 0
\end{bmatrix},\ Y_{n+i}^n=
\begin{bmatrix}
  0 & \cdots & 0\\
  \vdots & & \vdots \\
  0 & \cdots & 0\\
  \hline
  & I_n &
\end{bmatrix}
\in\matring_{n+i,n}(K),$$
both having $i$ zero rows of length $n$.
If $Z'$ is some matrix, then we call the matrix
$$ Z=
\left[\begin{array}{ccc|ccc}
  & & & 1 & & \\
    & Z' & & & \ddots & \\
  \hline
  & & & 1 & & 1\\
  & 0 & & & \ddots & \\
  & & & & & 1
\end{array}\right]$$ with entries $1$ on two diagonals each of length
$m \geq 0$ the $m$-th \emph{enlargement\/} of $Z'$.

A typical example of a series of preprojective indecomposable
modules of rank $2$ over a canonical algebra of type $(p,2,2)$ is
the following. We fix $i$ and$j$ with $ 1\leq i < j \leq p$ and
consider the module

\vspace{0,3cm}
  $M_m^{(i,j)}$

\vspace{-0,8cm}

{\footnotesize
  $$\xy\xymatrixcolsep{1.0pc}\xymatrix{ { } & K^{m} \ar @{=}[r] &
    \dots \ar @{=}[r] & K^{m} \ar @{->}[r]^-{X^{m}_{m+1}} & K^{m+1} \ar
    @{=}[r] & \dots \ar @{=}[r] & K^{m+1} \ar @{->}[r]^-{X^{m+1}_{m+2}} &
    K^{2m+2} \ar @{=}[r] & \dots \ar @{=}[r] & K^{m+2} \ar
    @{=}[dr] &\\
    K^m \ar @{=}[ur] \ar @{->}[drrrrr]_-{Y^{m}_{m+1}} \ar
    @{->}[rrrrr]^-{Y^{m}_{m+1}} & &
    & & & K^{m+1} \ar @{->}[rrrrr]^-{Y^{m+1}_{m+2}} & & & & & K^{m+2}
    \\ & & & & &
    K^{m+1} \ar @{->}[rrrrru]_-{Z^{m+1}_{m+2}} & & & & & } \endxy
$$}

\noindent where $Z^{m+1}_{m+2} $ is the m-th enlargement of the  $2
\times 1$ matrix $ Z'=
\begin{bmatrix}
  1 \\
  1 \\
\end{bmatrix}$.
Here  the matrices  $ X_{m}^{m+1}$ and $ X_{m+1}^{m+2}$   are
 associated
 to the arrows $ \al_i : ( i-1) \ra i$
and  $ \al_{j} : ( j-1)  \ra j $ respectively.
It follows from \cite{kussin:meltzer:04} that each preprojective
indecomposable rank $2$ module for a canonical algebra of type
$(p,2,2)$ is isomorphic to a module of this form.

The indecomposable modules for other domestic canonical algebras are
defined in a similar way by enlargement of certain ``small'' matrices
which can be found in~\cite[Theorem~2]{kussin:meltzer:04}
and~\cite{komoda:meltzer}, respectively.

Now, let $A=K\Gamma$ be the path algebra of an extended Dynkin quiver.
As we have already mentioned in the introduction
 there is a canonical algebra
of domestic type $\La$ and a left  tilting  module $T$ such that
$\End(T)\simeq A^{op}$.
For general information about tilting theory we refer to
\cite{brenner:butler}
and \cite{happel:ringel}.

We consider the functor $F=\Hom(T,-): \Lambda-\mod \rightarrow
\mod-A^{op}   $, where
$\ \Lambda-\mod$ (respectively
 $\mod-A^{op} $)   is the category of finite-dimensional left $\Lambda$
modules.
(respectively finite-dimensional right
$A^{op}$ modules).
Obviously the last  can be identified with  $A-\mod$, thus with the
category  of representations of  $\Gamma$.

We recall that a homomorphism $ f:M \rightarrow M'$ of
$\Lambda$-modules is given by a set of linear maps $ f_i: M(i)
\rightarrow M'(i)$ such that for each arrow $\phi :i \rightarrow j$ of
the quiver for $\Lambda$ we have $ f_j M(\phi) = M'(\phi) f_i$. From
\cite[Lemma 4.2]{geigle:lenzing} we know that the linear maps
$M(\phi)$ for a indecomposable preprojective module $M$ are
monomorphisms. As a consequence a homomorphism between indecomposable
preprojective modules $M$ and $M'$ is uniquely determined by the map
$f(\infty)$ and we will always identify such a homomorphism with the
matrix for $f(\infty)$.

We write the tilting module $T$ as a direct sum of pairwise
non-isomorphic indecomposables $T_1,\dots, T_l$ and choose
generators for all non-zero vector spaces $\Hom(T_i,T_j)$.  If such
a homomorphism space is non-zero, it is $1$-dimensional and in this
case  a generator $f$ can be identified with a a single matrix
$S^i_j$ describing the linear map $f({\infty}) $.

Now, if $M$ is a preprojective left $\La$-module the corresponding
representation $N=F(M)$ of the extended Dynkin quiver $\Gamma$ can
be computed as follows.  The vector spaces for the vertices $i$ of
$\Gamma$ are given by $N(i)=\Hom(T_i,M)$. Furthermore, for each $i$
we choose a basis of $N(i)$. In case $\Hom(T_i,T_j) \neq 0$ there is
an arrow $\phi: i \ra j $ in $\Gamma^{op}$.  We again have that
$\Hom(T_i,M)$ (respectively $\Hom(T_j,M)$) can be identified with a
vector space of matrices describing the linear map in $\infty$. Then
the linear map $N(\phi) : N(j)=\Hom(T_j,M) \ra \Hom(T_i,M)=N(i)$ is
the multiplication with the matrix $S^i_j$ from the right hand side.
Consequently, by our choice of bases in the $N(i)$ we obtain the
matrices for $N$. The rank of a representation of $\Gamma$ is by
definition the rank of the corresponding $\Lambda$-module. Note that
in~\cite{kussin:meltzer:04} (respectively~\cite{komoda:meltzer}) the
$\La$-modules $M$ are constructed as members of a series of
indecomposable $\La$-modules $M_m$ ($m\geq 0$). The procedure just
described will be applied simultaneously to the whole series.

The theorem of Brenner and Butler~\cite{brenner:butler} implies that
if $M$ satisfies the condition $\Ext^1(T,M)=0$ and $M$ is
indecomposable then $F(M)$ is also indecomposable.  We will apply
the functor $F$ to indecomposable preprojective $\La$-modules
satisfying the condition above. In this way we obtain the matrices
of all indecomposable preprojective representations of the extended
Dynkin diagram $\Gamma$.

Two kinds of data are important for our construction: the explicit
knowledge of the tilting module and the explicit knowledge of the
indecomposable preprojective $\La$-modules, both given by vector
spaces and matrices.

We note that the same method can be applied also to determine explicitly the
indecomposable modules over a tame concealed algebra.

\section{The case  $\widetilde{D}_n$ }\label{Dn}

\begin{Number}

For the structure of the indecomposable representations of extended
Dynkin quivers of type $\widetilde{A}_n$ we refer to
\cite[chapter~11]{gabriel:roiter}. In this chapter we study the
representations of the extended Dynkin quiver of type
$\widetilde{D}_n$ where we fix the following orientation.

\vspace{0,3cm}
\hspace{1,3cm}
$\Gamma$
\vspace{-0,7cm}

  {\footnotesize 
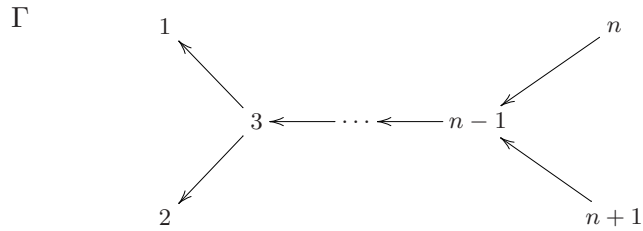
$$\xy\xymatrixcolsep{2.0pc}\xymatrix{
  1  &     &              &      & n   \ar @{->}[ld]  \\
    &    3 \ar @{->}[lu]  \ar @{->}[ld]  &    \dots  \ar @{->}[l]    &    n-1  \ar @{->}[l]    &           \\
  2  &     &                     &      & n+1  \ar @{->}[lu]   }
  \endxy
  $$}

  \vspace{0,2cm}

  The corresponding canonical algebra $\La$ is of type $ (n-2,2,2)$.
  Looking at the preprojective component of $\La$ we see that the
  following $T=\bigoplus_{k=1}^{n+1}T_k$ is a tilting module in
  $\La-\mod$ with $\End(T) \cong K\Gamma^{op}$.
  The indecomposable
  direct summands $T_k$ of $T$ are given as follows:

\newpage
For $3 \leq k \leq n-1$ we have

\vspace{0.2cm} \hspace{2cm}
$T_k:$
\vspace{-0,6cm}
 {\footnotesize $$\xy\xymatrixcolsep{1.0pc}\xymatrix{
 & 0    \ar @{->}[r]    &  \dots   \ar @{->}[r]     & 0 \ar @{->}[r]
 &   K      \ar @{->}[r]^-{=}
 &   \dots   \ar @{->}[r]^-{=}      & K    \ar @{->}[rd]^-{X^1_2 }     &     \\
0   \ar @{->}[rrrrd]    \ar @{->}[ru]   \ar @{->}[rrrr]   &      &
&        &
        K  \ar@{->}[rrr]^-{Y^1_2 }  &  &      &  K^2   \\
 &      &         &        &            K   \ar @{->}[rrru]_-{Z^1_2 }
 &   &         &     \\
  }
  \endxy
  $$}

\noindent
with $n-k-1$   entries $0$ and    $k-2$ entries $K$ in the first arm
(in particular for $T_{n-1}$ there is no $0$ in the first arm).
Moreover,
 {\footnotesize $$\xy\xymatrixcolsep{1.0pc}\xymatrix{
 & & 0    \ar @{->}[r]    &  \dots   \ar @{->}[r] & 0    \ar @{->}[rd]     &     \\
T_1: & 0   \ar @{->}[rrd]    \ar @{->}[ru]   \ar @{->}[rr]
&        &
        K  \ar@{->}[rr]^-{=}   &     &  K   \\
   &  &         &                 0   \ar @{->}[rru]
 &   &         &     \\
  }
  \endxy\quad \xy\xymatrixcolsep{1.0pc}\xymatrix{
 & & 0    \ar @{->}[r]    &  \dots   \ar @{->}[r] & 0    \ar @{->}[rd]     &     \\
T_2: & 0   \ar @{->}[rrd]    \ar @{->}[ru]   \ar @{->}[rr]
&        &
        0  \ar@{->}[rr]   &     &  K   \\
    &  &         &                 K   \ar @{->}[rru]^-{=}
 &   &         &     \\
  } \endxy$$}
and
{\footnotesize $$\xy\xymatrixcolsep{1.0pc}\xymatrix{
 & & 0    \ar @{->}[r]    &  \dots   \ar @{->}[r] & 0    \ar @{->}[rd]     &     \\
T_n : & 0   \ar @{->}[rrd]    \ar @{->}[ru]   \ar @{->}[rr]
&        &
        K  \ar@{->}[rr]^-{=}   &     &  K   \\
   &  &         &                 K   \ar @{->}[rru]^-{=}
 &   &         &     \\
  }
  \endxy\quad \xy\xymatrixcolsep{1.0pc}\xymatrix{
    & & K    \ar @{->}[r]^-{=}    &  \dots   \ar @{->}[r] & K    \ar
    @{->}[rd]^-{X_2^1}     &     \\
    T_{n+1}: & K \ar @{->}[rrd]^-{=} \ar @{->}[ru] \ar @{->}[rr]^-{=} & &
    K  \ar@{->}[rr]^-{Y_2^1}   &     &  K^2   \\
    & & & K \ar @{->}[rru]^-{Z_2^1}
    &   &         &     \\
  } \endxy$$}

Here $Z^1_2 $ denotes the matrix
$\begin{bmatrix}
  1 \\
  1 \\
\end{bmatrix}$.

The following picture indicates the quiver for the endomorphism ring
of $T$ and gives generators $S^i_j$ for the non-zero homomorphism
spaces
$\Hom_{\La}(T_i,T_j)$ which are always represented by the matrices  for the
linear maps $f_{\infty} : T_i(\infty) \ra T_j(\infty)   $.

\vspace{-0,3cm}
 $$\xy\xymatrixcolsep{2.3pc}\xymatrix{
 T_1  \ar @{->}[rd]^-{ \begin{bmatrix} 0  \\1  \end{bmatrix}   }   &     &     &       &    T_n  \\
   &   T_3 \ar @{->}[r]_-{=}      &  \ar @{->}[r]_-{=}   \dots  \ar @{->}[r]_-{=}   &
 T_{n-1}   \ar @{->}[ru]^-{ \begin{bmatrix} 0  & 1  \end{bmatrix}   }    \ar @{->}[rd]_-{=}      &      \\
 T_2  \ar @{->}[ru]_-{ \begin{bmatrix} 1  \\1  \end{bmatrix}   }      &     &     &       &   T_{n+1}   \\
}
  \endxy
$$
\end{Number}

\begin{Number}

  We start with the description of the indecomposable preprojective
  modules of rank $2$.  There are precisely $\binom{ n-2}{2}$ series
  of indecomposable $\La$-modules of rank $2$: $M_m^{(i,j)}$, $1 \leq i
  < j \leq n-2$. They have been described in the previous chapter.

We now fix $i$ and $j$ with $1 \leq i <  j \leq n-2$ and compute
simultaneously the representations   $ N_m^{(i,j)} =  F(  M_m^{(i,j)} )$, $
  m \in \NN$. We shortly write $ M_m^{(i,j)}=M$.

\vspace{0,2cm}
\underline{case (a)}
We assume  that $i\neq 1 $ and $ j \neq n-2$.

 (a1)
Computation of  $\Hom(T_1,M)$:

A homomorphism $f: T_1 \ra M$ is given by matrices $Q=(q_i)
 \in\matring_{m+1,1}(K)$
 and $S=(s_i)\in \matring_{m+2,1}(K)$ such that
$S=Y^{m+1}_{m+2} Q$.
We have already mentioned that for homomorphisms between
indecomposable preprojective modules
 a homomorphism  is uniquely
determined
by the matrix for the linear map of the point $\infty$, that is
  $S$.
The matrix equation yields that
 $s_1=0$  and $ s_l=q_{l-1}$ for $l=2,\dots, m+2$.
Therefore $\dim_K\Hom(T_1,M)=m+1$ and a basis is given by $(m+2) \times
1$-matrices
$w^{(1)}_{2}, w^{(1)}_{3}   \dots,w^{(1)}_{m+2}$, where $w^{(1)}_{i}$ is
the matrix  with  entries $s_i=1$ and
$s_j=0$ for $j \neq i$.

\vspace{0,2cm}
(a2)
Computation of  $\Hom(T_2,M)$:

A homomorphism $f: T_1 \ra M$ is given by matrices $R=(r_i)
 \in\matring_{m+1,1}(K)$
 and $S=(s_i)\in \matring_{m+2,1}(K)$ such that
$S=Z^{m+1}_{m+2} R$.

Now, the matrix equation and the shape of  $Z^{m+1}_{m+2}$ imply that
{\small
 $$S=\begin{bmatrix}
r_1+r_2  \\
 r_1+r_3 \\
r_2+r_4 \\
r_3+r_5 \\
  \vdots  \\
\end{bmatrix}.$$
}

\noindent
(We formally define $r_i=0$ for $i>m+1$ and $s_i=0$ for $i>m+2$.)
Observe that  then {\small
 $$r_1+r_2=(r_1+r_3)+(r_2+r_4)-(r_3+r_5)-(r_4+r_6)+(r_5+r_7)+(r_6+r_8)-\dots
 $$}
\noindent
 and consequently
{\small
$$ s_1=s_2+s_3-s_4-s_5+s_6+s_7-\dots .$$
}
\noindent
Therefore $\dim_K\Hom(T_2,M)=m+1$ and a basis is given by $(m+2) \times
1$-matrices
$w^{(2)}_{2}, w^{(2)}_{3}    \dots, w^{(2)}_{m+2}$, where
 $w^{(2)}_{i}$ is
the matrix  with  entries $s_i=1$, $s_j=0$ for $j \neq 1,  i$ and
$s_1$ is $1$ or $-1$, which is dependent on  the rest of $i$ modulo
$4$.

\vspace{0,2cm}
(a3)
Computation of  $\Hom(T_k,M)$ for $ k=3,\dots,n-j$:

A homomorphism $f: T_k \ra M$ is given by matrices $P=(p_i)
\in\matring_{m+2,1}(K)$, $Q=(q_i) \in\matring_{m+1,1}(K)$, $R=(r_i)
\in\matring_{m+1,1}(K)$ and $S=(s_{i,j})\in \matring_{m+2,2}(K)$
such that $SX^{1}_{2}=P$, $SY^{1}_{2} = Y^{m+1}_{m+2}  Q$ and
$SZ^{1}_{2}=Z^{m+1}_{m+2}R$.

The first equation yields no condition for the coefficients of $S$
whereas from the second equation we conclude that $s_{1,2}=0$.
The third condition shows as in the case (a2) that the coefficient
$s_{1,1}$ is a linear combination of the remaining coefficients.

Therefore $\dim_K\Hom(T_k,M)=m+2$ and a basis is given by $(m+2) \times
2$-matrices
$w^{(k)}_{2,1}, w^{(k)}_{3,1}, \dots, w^{(k)}_{m+2,1},
w^{(k)}_{2,2}, w^{(k)}_{3,2}, \dots, w^{(k)}_{m+2,2}$,
 where
 $w^{(k)}_{i,j}$ is
the matrix  with  entries $s_{i,j}=1$, $s_{u,v}=0$ for $(u,v) \neq
(i,j)$ and $(1,1)$
 whereas  $s_{1,1}$ is $1$ or $-1$, which is
dependent on the the rest of $i$ modulo $4$.

\vspace{0,2cm}
(a4)
Computation of  $\Hom(T_k,M)$ for $ k=n-j+1,\dots,n-i$:

In this case a homomorphism $f: T_1 \ra M$ is given by matrices
$P'=(p'_i) \in\matring_{m+1,1}(K)$, $P=(p_i) \in\matring_{m+2,1}(K)$,
$Q=(q_i) \in\matring_{m+1,1}(K)$, $R=(r_i) \in\matring_{m+1,1}(K)$ and
$S=(s_{i,j})\in \matring_{m+2,2}(K)$ such that $P=X^{m+1}_{m+2} P'$,
$SX^{1}_{2}=P$, $SY^{1}_{2} = Y^{m+1}_{m+2} Q$ and
$SZ^{1}_{2}=Z^{m+1}_{m+2}R$.

In this case the first equation yields $p_{m+2}=0$ and together with
the second equation we get $s_{m+2,1}=p_{m+2}=0$. As in the case
above  we conclude that $s_{1,2}=0$ and that $s_{1,1}$ is a linear
combination of the remaining coefficients.

Therefore $\dim_K\Hom(T_k,M)=m+1$ and a basis is given by $(m+2) \times
2$-matrices
$w^{(k)}_{2,1}, w^{(k)}_{3,1}, \dots,w^{(k)}_{m+1,1},
w^{(k)}_{2,2}, w^{(k)}_{3,2}, \dots,w^{(k)}_{m+2,2}$,
 where
 $w^{(k)}_{i,j}$ is
the matrix  with  entries $s_{i,j}=1$ and $s_{u,v}=0$ for $(u,v)
\neq (i,j)$ and $(1,1)$ whereas  $s_{1,1}$ is $1$ or $-1$.

\vspace{0,2cm}
(a5)
Computation of  $\Hom(T_k,M)$ for $ k=n-i+1,\dots,n-1$:

In this case a homomorphism $f: T_k \ra M$ is given by matrices
$P''=(p''_i) \in\matring_{m,1}(K)$,
$P'=(p'_i) \in\matring_{m+1,1}(K)$,
$P=(p_i) \in\matring_{m+2,1}(K)$,
$Q=(q_i) \in\matring_{m+1,1}(K)$,
$R=(r_i) \in\matring_{m+1,1}(K)$ and
$S=(s_{i,j})\in \matring_{m+2,2}(K)$ such that
$P'=X^{m}_{m+1} P''$
$P=X^{m+1}_{m+2} P'$
$SX^{1}_{2}=P$,
$SY^{1}_{2} = Y^{m+1}_{m+2}  Q$ and
$SZ^{1}_{2}=Z^{m+1}_{m+2}R$.

The first equation yields $p'_{m+1}=0$ and together with the second
equation we get  $p_{m+1}=p'_{m+1}=0$ and additionally  $p_{m+2}=0$.
Then from $SX^{1}_{2}=P$ we conclude $s_{m+1,1}=s_{m+2,1}=0$. As in
the  (a3) the  other equations imply that $s_{1,2}=0$ and that
$s_{1,1}$ is a linear combination of the remaining coefficients.

Therefore $\dim_K\Hom(T_k,M)=m$ and a basis is given by $(m+2) \times
2$-matrices
$w^{(k)}_{2,1}, w^{(k)}_{3,1}, \dots,w^{(k)}_{m,1},
w^{(k)}_{2,2}, w^{(k)}_{3,2}, \dots,w^{(k)}_{m+2,2}$,
 where
 $w^{(k)}_{i,j}$ is
the matrix  with  entries $s_{i,j}=1$ and $s_{u,v}=0$ for $(u,v)
\neq (i,j)$ and $(1,1)$
 whereas $s_{1,1}$ is $1$ or $-1$.

\vspace{0,2cm}
(a6)
Computation of  $\Hom(T_n,M)$:

A homomorphism $f: T_n \ra M$ is given by matrices
$Q=(q_i) \in\matring_{m+1,1}(K)$,
$R=(r_i) \in\matring_{m+1,1}(K)$ and
$S=(s_{i,j})\in \matring_{m+2,2}(K)$ such that
$S=Y^{m+1}_{m+2} Q$ and
$S=Y^{m+1}_{m+2} R$.

The first equation implies $ s_1=0$ and
 the second equation gives
{\small
$$
S=\begin{bmatrix}
0\\ s_2\\s_3\\s_4\\ \vdots
\end{bmatrix}
=
\begin{bmatrix}
r_1+r_2 \\r_1+r_3\\r_2+r_4\\r_3+r_5\\ \vdots
\end{bmatrix}
$$
}

Now
$$0=(r_1+r_2)=(r_1+r_3)+(r_2+r_4)-(r_3+r_5)-(r_4+r_6)+(r_5+r_7)+
\dots   \, . $$
\noindent
which  implies
$$0 =s_2+s_3-s_4-s_5+s_6+\dots$$
\noindent
and $s_2$ can be written as a linear combination of the remaining
coefficients.
(We again define formally $r_i=0$ for $i>m+1$
                      and $s_i=0$ for $i>m+2$).

                      Therefore $\dim_K\Hom(T_n,M)=m$ and a basis is
                      given by $(m+2) \times 1$-matrices $w^{(n)}_ 3 ,
                      w^{(n)}_4 , \dots,w^{(n)}_{m+2}$ where $w^{(n)}_
                      i$ is the matrix with entries $s_1 =0$, $s_i=1$,
                      $s_j=0$ for $j\neq 1,2, i$ and $s_2=1$ if $i
                      \equiv 0,1 \mod 4 $ and $s_2= - 1$ if $i \equiv
                      2,3 \mod 4 $.

\vspace{0,2cm}
(a7)
Computation of  $\Hom(T_{n+1},M)$:

A  homomorphism $f: T_{n+1} \ra M$ is given by matrices
$U=(u_{i})\in \matring_{m,1}(K)$
$P'=(p'_i) \in\matring_{m+1,1}(K)$,
$P=(p_i) \in\matring_{m+2,1}(K)$,
$Q'=(q'_i) \in\matring_{m+1,1}(K)$,
$Q=(q_i) \in\matring_{m+2,1}(K)$,
$R=(r_i) \in\matring_{m+1,1}(K)$ and
$S=(s_{i,j})\in \matring_{m+2,2}(K)$ such that
$P'=X^{m}_{m+1} U$,
$P=X^{m+1}_{m+2} P'$,
$SX^{1}_{2}=P$,
$Q'=Y^{m}_{m+1} U$,
$Q=Y^{m+1}_{m+2} Q'$,
$SY^{1}_{2} = Y^{m+1}_{m+2}  Q$,
$R=Y^{m}_{m+1} U$
 and
$SZ^{1}_{2}=Z^{m+1}_{m+2}R$.

It is easily calculated that the equations imply that $S$ is of the
form
$$S=\begin{bmatrix}
u_1  &0 \\
u_2  &0 \\
u_3  &u_1 \\
\vdots & \vdots \\
u_m & u_{m-2} \\
0 &  u_{m-1} \\
0 &  u_{m} \\
\end{bmatrix}
.$$

Therefore $\dim_K\Hom(T_{n+1},M)=m$ and a basis is given by $(m+2) \times
2$-matrices
$w^{(n+1)}_1  , w^{(n+1)}_2  , \dots,  w^{(n+1)}_m$  ,
where  $ w^{(n+1)}_i$ is the $(m+2)\times 2$-matrix  with  entries
$s_{i,1}=1$, $s_{i+2,2}=1$ and all other entries are zero.

\vspace{0,2cm} Now, in order to determine the matrices of the
representation $N=F(M)$ we have to describe the linear maps
$N(j)=\Hom(T_j,M) \ra \Hom(T_i,M) =N(i) $ in the given bases.  As
mentioned in Section~\ref{sec:tilting} this map is identified with the
multiplication of the matrix $S^j_i$ from the right.

In particular the map $N(n) \ra N(n-1)$ is given by the formula
{\small
$$\begin{bmatrix}
0 \\ s_2 \\s_3 \\ \vdots \\ s_{m+2}
\end{bmatrix}
\mapsto
\begin{bmatrix}
0 \\ s_2 \\s_3 \\ \vdots \\ s_{m+2}
\end{bmatrix}
\begin{bmatrix}
0 &1
\end{bmatrix}
=
\begin{bmatrix}
0 &0\\ 0 &  s_2 \\ 0 & s_3 \\ \vdots & \vdots \\ 0 &  s_{m+2}
\end{bmatrix}.
$$
}

Thus we obtain in the bases
$w_3^{(n)},w_4^{(n)} \dots, w^{(3)}_{m+2}$
 and
$w^{(n-1)}_{2,1}, w^{(n-1)}_{3,1},   \dots, w^{(n-1)}_{m,1}, \\
  w^{(n-1)}_{2,2},   w^{(n-1)}_{3,2},   \dots,
 w^{(n-1)}_{m+2,2}$
the following matrix

{\footnotesize
$$ C=
\begin{array}{  |   cccc  cccc  c |  }
\hline
  & & &    & &     & &     &  \\
  & & &    & &     & &     &  \\
  & & &    & &     & &     &   \\
  & & &    & &     & &     &  \\
  & & &    & &     & &     &  \\
  & & &    & &     & &     &   \\
\hline
 -1 &1 &1 &-1    &-1 & 1    &1 &  -1   & \dots  \\
  &1 & &    & &     & &     &  \\
  & &1 &    & &     & &     &   \\
  & & & 1   & &     & &     &  \\
  & & &    &1 &     & &     &  \\
  & & &    & & \ddots    & &     &   \\
\hline
\end{array}
$$
}

All the other matrices for $N$ are computed in the same way and we
obtain  in the case (a)  the following representation
$N^{(i,j)}_m=F(M^{(i,j)}_m)$

{\footnotesize $$\xy\xymatrixcolsep{1.0pc}\xymatrix{ K^{m+1} & & & & & &
    & & & & K^m
    \ar @{->}[dl]^-{C} \\
    { } & K^{2m+2} \ar @{->}[ul]_-{A} \ar @{->}[dl]^-{B} & \cdots \ar
    @{->}[l]_-{=} & K^{2m+2} \ar @{->}[l]_-{=} & K^{2m+1} \ar
    @{->}[l]_-{E} & \cdots \ar @{->}[l]_-{=} & K^{2m+1} \ar
    @{->}[l]_-{=} & K^{2m} \ar @{->}[l]_-{F}
    & \cdots \ar @{->}[l]_-{=} & K^{2m} \ar @{->}[l]_-{=} & \\
    K^{m+1} & & & & & & & & & & K^m \ar @{->}[ul]^-{D} } \endxy
$$}
with
\begin{gather}\label{eq:matrix-anfang}
A=\boxed{\begin{array}{c|c}
 \mathbf{0} &  \mathbf{I_{m+1} }
\end{array}}\quad B=\boxed{\begin{array}{c|c}
  \mathbf{I_{m+1}} &  \mathbf{I_{m+1} }
\end{array}}\quad C=\boxed{\begin{array}{c}
 \mathbf{0} \\
 \hline
 v \\
  \hline
  \mathbf{ I_{m} }
\end{array}}\quad D=\boxed{\begin{array}{c}
 \mathbf{0} \mid  \mathbf{I_{m-1}}\\
 \hline
 0\dots 0\\
 \hline
  \mathbf{I_m}
\end{array}}\\
E=\boxed{\begin{array}{c|c}
 \begin{array}{c} \mathbf{ I_{m}}\\ \hline 0\dots 0\end{array} & \\
 \hline
  &  \mathbf{I_{m+1}}
\end{array}}\quad F=\boxed{\begin{array}{c|c}
 \begin{array}{c}  \mathbf{I_{m-1}}\\ \hline 0\dots 0\end{array} & \\
 \hline
  & \mathbf{ I_{m+1}}
\end{array}}
\end{gather}
where in $C$ the vector $v$ is given by the first $m$ entries of the
periodic vector
$$(-1\,1\,1\,-1\mid -1\,1\,1\,-1\mid\ldots).$$

In the ''degenerated cases'' similar calculations as above lead to the
following results

\underline{ case (b)}: $i=1$ and $j < n-2$

{\footnotesize 
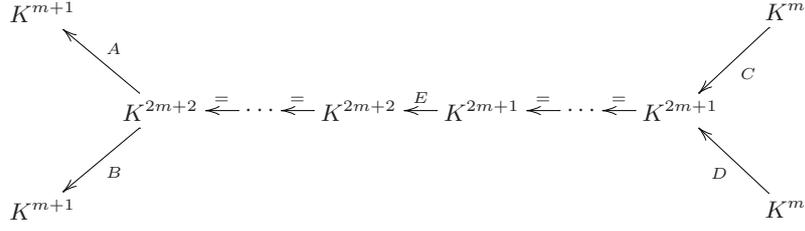
$$\xy\xymatrixcolsep{1.0pc}\xymatrix{ K^{m+1} & & & & & &
    & K^m
    \ar @{->}[dl]^-{C} \\
    { } & K^{2m+2} \ar @{->}[ul]_-{A} \ar @{->}[dl]^-{B} & \cdots \ar
    @{->}[l]_-{=} & K^{2m+2} \ar @{->}[l]_-{=} & K^{2m+1} \ar
    @{->}[l]_-{E} & \cdots \ar @{->}[l]_-{=} & K^{2m+1} \ar
    @{->}[l]_-{=} & \\
    K^{m+1} & & & & & & &  K^m \ar @{->}[ul]^-{D} } \endxy
$$}

with matrices $A$, $B$ and $E$ as given in case~(a) and with
\begin{equation}
  C=\boxed{\begin{array}{c}
 \mathbf{0} \\
 \hline
 v \\
  \hline
   \mathbf{I_{m} }
\end{array}}\quad   D=\boxed{\begin{array}{c}
 \mathbf{0} \mid  \mathbf{I_{m-1}}\\
 \hline
 0\dots 0\\
 \hline
 0\dots 0\\
 \hline
  \mathbf{I_m}
\end{array}}
\end{equation}
with vector $v$ as in case~(a).

\underline{ case (c)}: $i \neq 1$ and $j =  n-2$
{\footnotesize 
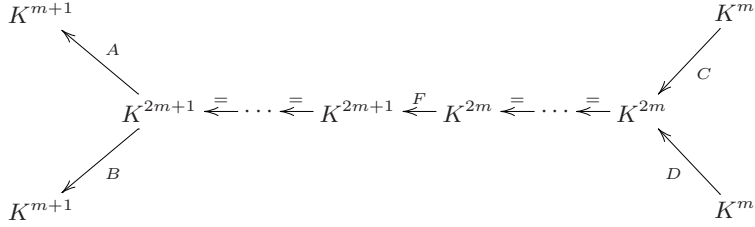
$$\xy\xymatrixcolsep{1.0pc}\xymatrix{ K^{m+1} & & &
    & & & & K^m
    \ar @{->}[dl]^-{C} \\
    { } & K^{2m+1} \ar @{->}[ul]_-{A} \ar @{->}[dl]^-{B} & \cdots \ar
    @{->}[l]_-{=} & K^{2m+1} \ar @{->}[l]_-{=} & K^{2m} \ar @{->}[l]_-{F}
    & \cdots \ar @{->}[l]_-{=} & K^{2m} \ar @{->}[l]_-{=} & \\
    K^{m+1} & & & & & & &  K^m \ar @{->}[ul]^-{D} } \endxy
$$}
with matrices $C$, $D$ and $F$ given as in case~(a), and with
\begin{equation}
  A=\boxed{\begin{array}{c|c}
 \mathbf{0} &  \mathbf{I_{m+1} }
\end{array}}\quad B=\boxed{\begin{array}{c|c}
  \begin{array}{c}
     \mathbf{I_{m}}\\
    \hline
    0\dots 0
  \end{array} &  \mathbf{I_{m+1} }
\end{array}}
\end{equation}

\underline{ case (d)}: $i =1$ and $j =  n-2$
{\footnotesize 
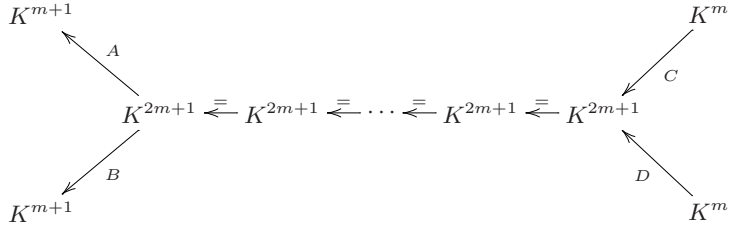
$$\xy\xymatrixcolsep{1.0pc}\xymatrix{ K^{m+1} & & & & &
    & K^m
    \ar @{->}[dl]^-{C} \\
    { } & K^{2m+1} \ar @{->}[ul]_-{A} \ar @{->}[dl]^-{B} & K^{2m+1}
    \ar @{->}[l]_-{=} & \cdots \ar @{->}[l]_-{=} & K^{2m+1} \ar
    @{->}[l]_-{=} & K^{2m+1} \ar @{->}[l]_-{=} & \\
    K^{m+1} & & & & & & K^m \ar @{->}[ul]^-{D} } \endxy
$$}
with matrices $A$, $B$  as in case~(c) and  and $C$, $D$ as in
case~(b).
\end{Number}

\begin{Number}
  In order to describe the representations of $Q$ of rank one we use
  the symmetry of the quiver. For this reason we study first by
  investigating the dimension vectors to which of the preprojective
  $\La$-modules we have to apply the functor $F$
  in order to get, up to symmetry, all preprojective indecomposable
  representations of $\Gamma$ of rank $1$.
  Recall that a
  tilting module $T$  induces an isomorphism of the corresponding
  Grothendieck groups such that the dimension vectors of the
  indecomposable direct summands $T_i$ are mapped to the dimension
  vectors of the left indecomposable projective modules $P_i$ over
  $A=\End(T)$ (see~\cite[3.2]{happel:ringel}). In our case
   this isomorphism is given by
$$f: \ZZ^{n+1}\cong K_0(\La) \ra K_0(A) \cong \ZZ^{n+1}$$
\vspace{-0,4cm}
{\footnotesize
$$
\begin{bmatrix}
   & a_1 & a_2 & \dots      &a_{n-3} &  \\
a_0&    &     &   a_{1'}    &        & a_c  \\
   &    &     &     a_{1''} &        &
\end{bmatrix}
\mapsto
\begin{bmatrix}
 a_{1'}&   &  &   &       &   x \\
       & a_{n-3}+x & \dots & a_2+x & a_1+x &  \\
a_{1''}&   &  &   &      &    x \\
\end{bmatrix}
$$
}
\noindent
where $x= a_{1'}+a_{1''}-a_c$.

There are $4(n-2)$ series of rank $1$-modules over $\La$ and the
following table shows how their dimension vectors are mapped under
the isomorphism $f$

{\footnotesize
$$
\begin{array}{ |c|c|c|c| }
\hline
type  & notation   &   \underline{\dim} (M)   & f(  \underline{\dim} (M)  )  \\
\hline
\hline
 1 & _{(1)}M^{(i)}_m    &
\begin{bmatrix}
   & \dots & m  & m+1  & \dots  &  \\
m &    &     &   &   m     &  m+1 \\
   &    &     &     &    m+1    &
\end{bmatrix}
  &
   \begin{bmatrix}
 m  &   &  &   &          m \\
       & 2m+1 & \dots & 2m  &  \\
m+1   &   &  &   &         m \\
\end{bmatrix}    \\
\hline
 2 & _{(2)}M^{(i)}_m    &
\begin{bmatrix}
   & \dots & m  & m+1  & \dots  &  \\
m &    &     &   &   m+1     &  m+1 \\
   &    &     &     &    m+1    &
\end{bmatrix}
  &
   \begin{bmatrix}
m+1  &   &  &   &          m+1 \\
       & 2m+2 & \dots & 2m+1  &  \\
m+1   &   &  &   &          m \\
\end{bmatrix}    \\
\hline
 3 & _{(3)}M^{(i)}_m    &
\begin{bmatrix}
   & \dots & m  & m+1  & \dots  &  \\
m &    &     &   &   m     &  m+1 \\
   &    &     &     &    m    &
\end{bmatrix}
  &
   \begin{bmatrix}
 m  &   &  &   &         m-1 \\
       & 2m & \dots & 2m-1  &  \\
m   &   &  &   &          m \\
\end{bmatrix}    \\
\hline
 4 & _{(4)}M^{(i)}_m    &
\begin{bmatrix}
   & \dots & m  & m+1  & \dots  &  \\
m &    &     &   &   m+1     &  m+1 \\
   &    &     &     &    m    &
\end{bmatrix}
  &
   \begin{bmatrix}
m+1  &   &  &   &           m \\
       & 2m+1 & \dots & 2m  &  \\
m   &   &  &   &           m \\
\end{bmatrix}    \\
\hline
\end{array}
$$
}

Here the growth of the dimension from $m$ to $m+1$
  in the first arm  for $\underline{\dim}( M)    $
 is realized for the
arrow $\alpha_i : (i-1) \ra i $ which implies that the growth of the
dimension  in the middle part  of  $f(  \underline{\dim} (  M)  ) $
 is realized
for the arrow $ (n-i-1) \ra (n-i)$. Moreover, for all types we have
$ 1 \leq i \leq n-2$ which
 means that in the particular cases $i=1$ or $i=n-2$
there is no growth of the dimension in the middle part  of
   $ f(  \underline{\dim} (  M) )  $.

Because each preprojective representation is exceptional and
each exceptional $A$-module is uniquely determined,
 up to isomorphism, by its dimension vector,
using the symmetry of the quiver $\widetilde{D}_n$
 it is sufficient
to describe the modules
  $  _{(1)}N^{(i)}_m =    F(   _{(1)}M^{(i)}_m )  $ and
  $    _{(2)}N^{(i)}_m =     F( _{(2)}M^{(i)}_m  )  $,
both for  $ 1 \leq i \leq n-2$ and $m \in \NN$.
The  calculations are done in the same way as in the case of rank
 $2$-modules and lead to the following results:

 $  _{(1)}N^{(i)}_m $, $ i \neq 1, n-2$
{\footnotesize 
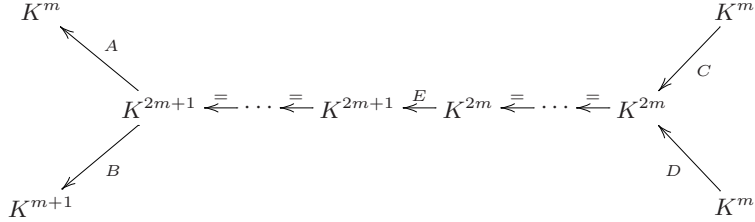
$$\xy\xymatrixcolsep{1.0pc}\xymatrix{ K^{m} & & &
    & & & & K^m
    \ar @{->}[dl]^-{C} \\
    { } & K^{2m+1} \ar @{->}[ul]_-{A} \ar @{->}[dl]^-{B} & \cdots \ar
    @{->}[l]_-{=} & K^{2m+1} \ar @{->}[l]_-{=} & K^{2m} \ar @{->}[l]_-{E}
    & \cdots \ar @{->}[l]_-{=} & K^{2m} \ar @{->}[l]_-{=} & \\
    K^{m+1} & & & & & & &  K^m \ar @{->}[ul]^-{D} } \endxy
$$}
with matrices
\begin{gather}
  A=\boxed{\begin{array}{c}
 \mathbf{0} \mid  \mathbf{I_{m}}
\end{array}}\quad B=\boxed{
\begin{array}{c|c}
  \mathbf{ I_{m+1} }&
  \begin{array}{c}
    0\dots 0\\
    \hline
   \mathbf{  I_m}
  \end{array}
\end{array}}\\
C=\boxed{\begin{array}{c}
 \mathbf{0} \\
 \hline
  \mathbf{ I_{m} }
\end{array}}\quad D=\boxed{\begin{array}{c}
  \mathbf{I_m }\\
 \hline
   \mathbf{I_{m} }
\end{array}}\quad E=\boxed{\begin{array}{c|c}
 \begin{array}{c} \mathbf{ I_{m}}
 \\ \hline 0\dots 0\end{array} & \\
 \hline
  &  \mathbf{I_{m}}
\end{array}}
\end{gather}

$ _{(1)}N^{(1)}_{m} $ {\footnotesize
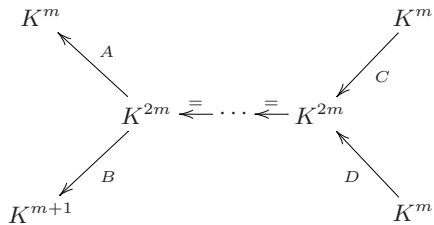
  $$\xy\xymatrixcolsep{1.0pc}\xymatrix{ K^{m} & & & & K^m
    \ar @{->}[dl]^-{C} \\
    { } & K^{2m} \ar @{->}[ul]_-{A} \ar @{->}[dl]^-{B} & \cdots \ar
    @{->}[l]_-{=}
    & K^{2m} \ar @{->}[l]_-{=} & \\
    K^{m+1}  & & & & K^m \ar @{->}[ul]^-{D} } \endxy
$$}
with
\begin{gather}
  A=\boxed{\mathbf{0}\mid  \mathbf{I_m}}
  \quad B=\boxed{
\begin{array}{c|c}
  \begin{array}{c}
    \mathbf{ I_{m} }\\
    {  }\\
    \hline
    0\dots 0
  \end{array} &
  \begin{array}{c}
    0\dots 0\\
    \hline
    {  }\\
    \mathbf{ I_m}
  \end{array}
\end{array}}
\end{gather}

 $  _{(1)}N^{(n-2)}_{m} $
{\footnotesize
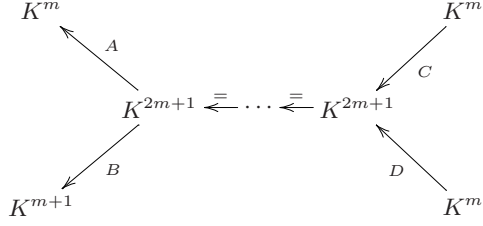
  $$\xy\xymatrixcolsep{1.0pc}\xymatrix{ K^{m} & & & & K^m
    \ar @{->}[dl]^-{C} \\
    { } & K^{2m+1} \ar @{->}[ul]_-{A} \ar @{->}[dl]^-{B} & \cdots \ar
    @{->}[l]_-{=}
    & K^{2m+1} \ar @{->}[l]_-{=} & \\
    K^{m+1}  & & & & K^m \ar @{->}[ul]^-{D} } \endxy
$$}
with
\begin{gather}
  C=\boxed{
  \begin{array}{c}
    \mathbf{0}\\
    \hline
     \mathbf{I_m}
  \end{array}}\quad D=\boxed{
  \begin{array}{c}
     \mathbf{I_m}\\
    \hline
    0\dots 0\\
    \hline
     \mathbf{I_m}
  \end{array}}
\end{gather}

 $  _{(2)}N^{(i)}_m $, $ i \neq 1, n-2$
{\footnotesize $$\xy\xymatrixcolsep{1.0pc}\xymatrix{ K^{m+1} & & &
    & & & & K^{m+1}
    \ar @{->}[dl]^-{C} \\
    { } & K^{2m+2} \ar @{->}[ul]_-{A} \ar @{->}[dl]^-{B} & \cdots \ar
    @{->}[l]_-{=} & K^{2m+2} \ar @{->}[l]_-{=} & K^{2m+1} \ar @{->}[l]_-{E}
    & \cdots \ar @{->}[l]_-{=} & K^{2m+1} \ar @{->}[l]_-{=} & \\
    K^{m+1} & & & & & & &  K^{m} \ar @{->}[ul]^-{D} } \endxy
$$}
with
\begin{gather}
  A=\boxed{\mathbf{0}\mid  \mathbf{I_{m+1}}}
  \quad B=\boxed{
\begin{array}{c|c}
  \begin{array}{c}
  \mathbf{   I_{m}} \\
    \hline
    0\dots 0
  \end{array} &  \mathbf{I_{m+1}}
\end{array}}\\
C=\boxed{\begin{array}{c}
  \mathbf{0} \\
 \hline
  \mathbf{ I_{m+1} }
\end{array}}\quad D=\boxed{
  \begin{array}{c}
     \mathbf{I_m}\\
    \hline
    0\dots 0\\
    \hline
     \mathbf{I_m}
  \end{array}}\quad E=\boxed{\begin{array}{c|c}
 \begin{array}{c}  \mathbf{I_{m}}
 \\ \hline 0\dots 0\end{array} & \mathbf{0}\\
 \hline
  & \mathbf{ I_{m+1}}
\end{array}}
\end{gather}

 $  _{(2)}N^{(1)}_{m} $
{\footnotesize
  $$\xy\xymatrixcolsep{1.0pc}\xymatrix{ K^{m+1} & & & & K^{m+1}
    \ar @{->}[dl]^-{C} \\
    { } & K^{2m+1} \ar @{->}[ul]_-{A} \ar @{->}[dl]^-{B} & \cdots \ar
    @{->}[l]_-{=}
    & K^{2m+1} \ar @{->}[l]_-{=} & \\
    K^{m+1}  & & & & K^m \ar @{->}[ul]^-{D} } \endxy
$$}
with \begin{gather}
 A=\boxed{\mathbf{0}\mid  \mathbf{I_{m+1}}}
\quad B=\boxed{
\begin{array}{c|c}
  \begin{array}{c}
     \mathbf{I_{m}} \\
    \hline
    0\dots 0
  \end{array} &  \mathbf{I_{m+1}}
\end{array}}
\end{gather}

 $  _{(2)}N^{(n-2)}_{m} $
{\footnotesize
  $$\xy\xymatrixcolsep{1.0pc}\xymatrix{ K^{m+1} & & & & K^{m+1}
    \ar @{->}[dl]^-{C} \\
    { } & K^{2m+2} \ar @{->}[ul]_-{A} \ar @{->}[dl]^-{B} & \cdots \ar
    @{->}[l]_-{=}
    & K^{2m+2} \ar @{->}[l]_-{=} & \\
    K^{m+1}  & & & & K^m \ar @{->}[ul]^-{D} } \endxy
$$}
with \begin{gather}\label{eq:matrix-ende}
  C=\boxed{
  \begin{array}{c}
    \mathbf{0}\\
    \hline
    \mathbf{ I_{m+1}}
  \end{array}}\quad D=\boxed{
  \begin{array}{c}
     \mathbf{I_m}\\
    \hline
    0\dots 0\\
    0\dots 0\\
    \hline
     \mathbf{I_m}
  \end{array}}
\end{gather}

\end{Number}

\begin{Thm}
  The representations given by the matrices described
  in~\eqref{eq:matrix-anfang}--~\eqref{eq:matrix-ende} above form a
  complete list of nonisomorphic preprojective
  indecomposable representations for the quiver
  $\widetilde{D}_n$ with the chosen orientation.
\end{Thm}

\section{The case  $\widetilde{E}_6$ }

We consider a quiver $\Gamma$  of type  $\widetilde{E}_6$ with subspace
orientation

  {\footnotesize 
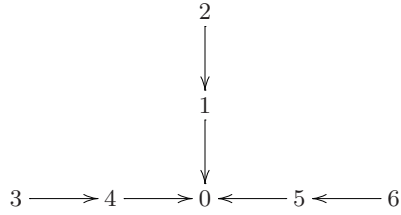
$$\xy\xymatrixcolsep{2.1pc}\xymatrix{
  &     &        2  \ar @{->}[d]    &      &  \\
         &      &  1   \ar @{->}[d]   &    &     \\
  3   \ar @{->}[r] &  4  \ar @{->}[r]  &      0       &  \ar
   @{->}[l]    5    &  \ar @{->}[l]   6  }
  \endxy
  $$}
In this case the corresponding domestic canonical algebra $\Lambda$ is
of type $(3,3,2)$ and a tilting module $T= \oplus_{i=0}^6 T_i$ \, such
that $\End (T) \simeq (K \Gamma)^{op}$ is given by

 {\footnotesize $$\xy\xymatrixcolsep{1.3pc}\xymatrix{
  T_0   &  K  \ar @{->}[r]^-{X^{1}_{2}}  & K^2   \ar @{->}[rd]^-{X^{2}_{3}}    &      \\
   0   \ar @{->}[r] \ar @{->}[ru]   \ar @{->}[rrd]   &  K  \ar @{->}[r]^-{Y^{1}_{2}}    & K^2  \ar @{->}[r]^-{Y^{2}_{3}}   & K^3    \\
     &      &  K    \ar @{->}[ru]_-{Z^{1}_{3}} &      \\
 }
  \endxy \quad
 \xy\xymatrixcolsep{1.3pc}\xymatrix{
  T_1   &  0   \ar @{->}[r]   & K  \ar @{->}[rd]^-{X^{1}_{2}}    &      \\
   0   \ar @{->}[r] \ar @{->}[ru]   \ar @{->}[rrd]   &
  K  \ar @{->}[r]^-{Y^{1}_{2}}    & K^2  \ar @{->}[r]^-{=}   & K^2   \\
     &      &  K    \ar @{->}[ru]_-{Z^{1}_{2}} &      \\
 }
  \endxy
\quad
 \xy\xymatrixcolsep{1.3pc}\xymatrix{
  T_2   &  0   \ar @{->}[r]   & K  \ar @{->}[rd]^-{=}    &      \\
   0   \ar @{->}[r] \ar @{->}[ru]   \ar @{->}[rrd]   &
  K  \ar @{->}[r]^-{=}    & K  \ar @{->}[r]^-{=}   & K  \\
     &      &  0    \ar @{->}[ru] &      \\
 }
  \endxy
  $$}

{\footnotesize $$\xy\xymatrixcolsep{1.3pc}\xymatrix{
 T_3   &  K   \ar @{->}[r]^-{X^{1}_{2}}   & K^2  \ar @{->}[rd]^-{=}    &      \\
   0   \ar @{->}[r] \ar @{->}[ru]   \ar @{->}[rrd]   &
  0  \ar @{->}[r]   & K  \ar @{->}[r]^-{Y^{1}_{2}}   & K^2   \\
     &      &  K    \ar @{->}[ru]_-{Z^{1}_{2}} &      \\
 }
  \endxy
\quad
 \xy\xymatrixcolsep{1.3pc}\xymatrix{
  T_4   &  K   \ar @{->}[r]^-{=}    & K  \ar @{->}[rd]^-{=}    &      \\
   0   \ar @{->}[r] \ar @{->}[ru]   \ar @{->}[rrd]   &
  0  \ar @{->}[r]  & K  \ar @{->}[r]^-{=}   & K  \\
     &      &  0    \ar @{->}[ru] &      \\
 }
  \endxy
  $$}

{\footnotesize $$\xy\xymatrixcolsep{1.3pc}\xymatrix{
 T_5  &  K   \ar @{->}[r]^-{=}   & K \ar @{->}[rd]^-{X^{1}_{2}}    &      \\
   0   \ar @{->}[r] \ar @{->}[ru]   \ar @{->}[rrd]   &
  K  \ar @{->}[r]^-{=}    & K  \ar @{->}[r]^-{Y^{1}_{2}}   & K^2   \\
     &      &  K    \ar @{->}[ru]_-{Z^{1}_{2}} &      \\
 }
  \endxy
\quad
 \xy\xymatrixcolsep{1.3pc}\xymatrix{
  T_6   &  K   \ar @{->}[r]^-{=}   & K \ar @{->}[rd]^-{X^{1}_{2}}    &      \\
   K   \ar @{->}[r]^-{=}  \ar @{->}[ru]^-{=}    \ar @{->}[rrd]^-{=}   &
  K  \ar @{->}[r]^-{=}    & K  \ar @{->}[r]^-{Y^{1}_{2}}   & K^2   \\
     &      &  K    \ar @{->}[ru]_-{Z^{1}_{2}} &      \\
 }
  \endxy
  $$}

The matrices $X^1_2$, $X^2_3$, $Y^1_2$, $Y^2_3$
are defined in the previous chapter and  $Z^1_3$ denotes the matrix
$ \begin{bmatrix}
 1\\ 1\\ 1
  \end{bmatrix}$.

 As in the
previous chapter we illustrate the endomorphism ring of $T$ and give
generators $S^i_j$ for the non-zero homomorphism spaces
$\Hom_{\La}(T_i,T_j)$.

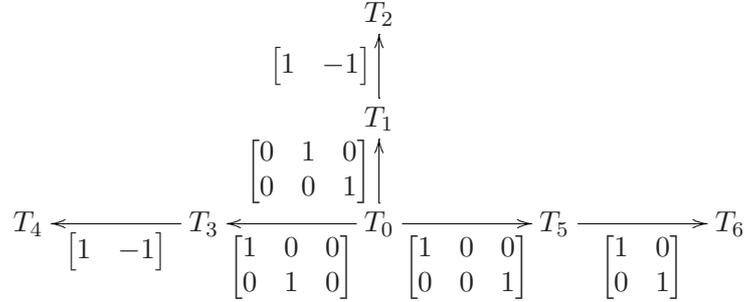
   $$\xy\xymatrixcolsep{4.pc}\xymatrix{
  &     &        T_2      &      &  \\
         &      &  T_1   \ar @{->}[u]^-{ \begin{bmatrix} 1 & -1  \end{bmatrix}   }    &    &     \\
  T_4    &  T_3  \ar @{->}[l]^-{ \begin{bmatrix} 1 & -1  \end{bmatrix}   }    &    \ar @{->}[l]^-{
  \begin{bmatrix} 1 & 0 & 0 \\ 0 & 1 &0 \end{bmatrix}   }   T_0   \ar @{->}[u]^-{
  \begin{bmatrix} 0 & 1 & 0 \\ 0 & 0 &1 \end{bmatrix}   }    \ar @{->}[r]_-{
  \begin{bmatrix} 1 & 0 & 0 \\ 0 & 0 &1 \end{bmatrix}   }     &  \ar
   @{->}[r]_-{
  \begin{bmatrix}  1 & 0 \\  0 &1 \end{bmatrix}   }    T_5    &    T_6  }
  \endxy
  $$

  In this chapter we determine matrices for the indecomposable
  preprojective representations corresponding to the indecomposable
  $\Lambda$-modules of rank $3$. There are 2 series of indecomposable
  preprojective rank $3$-modules over $\La$. The first one is
  described as follows: (see \cite{kussin:meltzer:04}.)

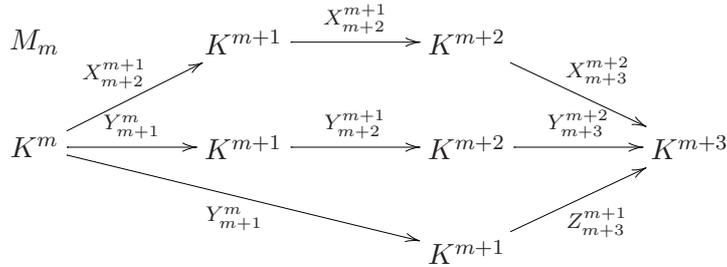
 $$
 \xy\xymatrixcolsep{4pc}\xymatrix{
  M_m   &  K^{m+1}   \ar @{->}[r]^-{X^{m+1}_{m+2}}       & K^{m+2} \ar @{->}[rd]^-{X^{m+2}_{m+3}}    &      \\
   K^m   \ar @{->}[r]^-{Y^{m}_{m+1}}       \ar@{->}[ru]^-{X^{m+1}_{m+2}}          \ar@{->}[rrd]_-{Y^{m}_{m+1}}      &
  K^{m+1}  \ar @{->}[r]^-{Y^{m+1}_{m+2}}  & K^{m+2}  \ar @{->}[r]^-{Y^{m+2}_{m+3}}   & K^{m+3}  \\
     &      &  K^{m+1}   \ar @{->}[ru]_-{Z^{m+1}_{m+3}} &      \\
 }
  \endxy
  $$
\noindent
where $  Z^{m+1}_{m+3}  $ is the $m$-th enlargement
 of the  $3 \times 1$ matrix $Z^1_3$.

As in chapter \ref{Dn} we apply the functor $F=\Hom_{\La}(T,-) $ to
the  modules $M_m$. We shortly write  $M_m=M$ and calculate first the
the dimensions and suitable bases of the vector spaces $N(i)$ of the
representation $ N=  \Hom_{\La}(T,M)$

\vspace{0,2cm}
 (a0)
Computation of  $\Hom(T_0,M)$:

 A homomorphism $f: T_0 \ra M$ is given by matrices
$P'=(p'_i) \in\matring_{m+1,1}(K)$,
$P =(p_i) \in\matring_{m+2,2}(K)$,
$Q'=(q'_i) \in\matring_{m+1,1}(K)$,
$Q =(q_{i,j}) \in\matring_{m+2,2}(K)$,
$R=(r_i) \in\matring_{m+1,1}(K)$ and
$S=(s_{i,j})\in \matring_{m+3,3}(K)$
such that
$P  X^{1}_{2}  =   X^{m+1}_{m+2}  P'$
$S  X^{2}_{3}=  X^{m+2}_{m+3}    P$,
$Q  Y^{1}_{2}  =   Y^{m+1}_{m+2}  Q'$
$S  Y^{2}_{3}=  Y^{m+2}_{m+3}    Q$,
 and
$SZ^{1}_{3}=Z^{m+1}_{m+3}R$.

It is easy to verify that the first four conditions imply that
$s_{1,2}=0$, $s_{1,3}=0$, $s_{2,3}=0$, $s_{m+2,1,1}=0$, $s_{m+3,1}=0$
and $s_{m+3,2}=0$, Moreover, the last equation gives that two of the
other entries of S (say $s_{1,1}$ and $s_{2,1}$) are linearly
dependent of the remaining $s_{i,j}$.  This is a consequence of the
structure of the matrix $Z^{m+1}_{m+3}$ and will appear in
calculations for all $N(i)$. However in order to determine the
matrices of the representation $N$ we need the precise description of
this linear dependence only for the vertex $1$ for which we give the
computation in detail. From this it will follow that also in the case
considered here we have the two linear dependent expressions mentioned
above. As a consequence we obtain that $\dim_K\Hom(T_0,M)=3m+1$ and a
basis is given by $(m+3) \times 3$-matrices $w^{(0)}_{3,1},
w^{(0)}_{4,1}, \dots,w^{(0)}_{m+1,1}, w^{(0)}_{2,2}, w^{(0)}_{3,2},
\dots,w^{(0)}_{m+2,2}, w^{(0)}_{3,3}, w^{(0)}_{4,3},
\dots,w^{(0)}_{m+3,3}, $ where $w^{(0)}_{i,j}$ is the matrix with
entries $s_{i,j}= 1$, with some possibly non-zero entries $s_{1,1}$,
$s_{2,1}$ and all the other entries are $0$.

\vspace{0,2cm}
 (a1)
Computation of  $\Hom(T_1,M)$:

 A homomorphism $f: T_1 \ra M$ is given by matrices
$P =(p_i) \in\matring_{m+2,1}(K)$,
$Q'=(q'_i) \in\matring_{m+1,1}(K)$,
$Q =(q_{i,j}) \in\matring_{m+2,2}(K)$,
$R=(r_i) \in\matring_{m+1,1}(K)$ and
$S=(s_{i,j})\in \matring_{m+3,2}(K)$
such that
$S  X^{1}_{2}=  X^{m+2}_{m+3}    P$,
$Q  Y^{1}_{2}  =   Y^{m+1}_{m+2}  Q'$
$S           =  Y^{m+2}_{m+3}    Q$,
 and
$SZ^{1}_{2}=Z^{m+1}_{m+3}R$.

 From the first conditions we conclude that   $s_{1,1}=0$,
 $s_{1,2}=0$,  $s_{2,2}=0$ and  $s_{m+3,1}=0$,
 whereas the last equation yields

{\small
$$ (1) \quad
\begin{bmatrix}
0  \\
s_{2,1} \\
s_{3,1}+ s_{3,2}  \\
s_{4,1}+ s_{4,2} \\
  \vdots  \\
\end{bmatrix}
=
\begin{bmatrix}
r_1+r_2  \\
r_1+r_3 \\
r_1+r_4 \\
r_2+r_5 \\
  \vdots  \\
\end{bmatrix}.
$$
}

We have the two identities
{\small
   $$ (2) \quad  r_1+r_2 =
  \sum_{k=0}^{\infty}  r_{6k+1}+r_{6k+4}
+ \sum_{k=0}^{\infty}  r_{6k+2}+r_{6k+5}
-  \sum_{k=0}^{\infty} r_{6k+4}+r_{6k+7}
-  \sum_{k=0}^{\infty} r_{6k+5}+r_{6k+8}
$$
}

{\small
 $$ (3) \quad   r_1+r_3 =
  \sum_{k=0}^{\infty}  r_{6k+1}+r_{6k+4}
+ \sum_{k=0}^{\infty}  r_{6k+3}+r_{6k+6}
-  \sum_{k=0}^{\infty} r_{6k+4}+r_{6k+7}
-  \sum_{k=0}^{\infty} r_{6k+6}+r_{6k+9}
$$
}

\noindent
(we formally  define  $ r_j=0$ for $j > m+1$).
According to $(1)$ and $(3)$  we get

{\small
 $$ (4) \quad   s_{2,1} =
  \sum_{k=0}^{\infty}  s_{6k+3,1}+s_{6k+3,2}
+ \sum_{k=0}^{\infty}  s_{6k+5,1}+s_{6k+5,2}
-  \sum_{k=0}^{\infty} s_{6k+6,1}+s_{6k+6,2}
-  \sum_{k=0}^{\infty} s_{6k+8,1}+s_{6k+8,2}
$$
}

\noindent (again we  formally  define  $ s_{i,j}=0$ for  $i>m+3$).
Further  it follows from $(1)$ and $(2)$ that one more coefficient
is dependent of the remaining, for instance we can write

\begin{eqnarray*}
    \label{eq:s31}
   (5) \quad  s_{3,1} & = & - s_{3,2}
- \sum_{k=1}^{\infty}  s_{6k+3,1}+s_{6k+3,2}
- \sum_{k=0}^{\infty}  s_{6k+4,1}+s_{6k+4,2} + \\
& &
+  \sum_{k=0}^{\infty} s_{6k+6,1}+s_{6k+6,2}
+  \sum_{k=0}^{\infty} s_{6k+7,1}+s_{6k+7,2}.
\end{eqnarray*}

It follows that  $\dim_K\Hom(T_1,M)=2m$ and a basis
$w^{(1)}_{4,1}, w^{(1)}_{5,1}   \dots,w^{(1)}_{m+2,1},
 w^{(1)}_{3,2},\\  w^{(1)}_{4,2}   \dots,w^{(1)}_{m+3,1}$
can is given in the following way: $w^{(1)}_{i,j}$  is the matrix
with entries $s_{i,j}=1$ and $s_{k,l}=0$ for $(k,l) \neq
(i,j),(2,1), (3,1)  $ and the entries $s_{3,1}$ and  $s_{2,1}$  have
to be computed using the formulas  $ (5)$ and  $ (4)$ (note that
first one has to calculate  $s_{3,1}$ because this coefficient
appears in  $ (4)$). The following table gives these coefficients
for the  vectors $w^{(1)}_{i,j}$

{\footnotesize
$$
\begin{array}{ |c|    c|c|c| c|c|c|    c ||    c|  c|c|c|    c|c|c|   c |   }
\hline
 vector &  4,1 & 5,1  & 6,1  &  7,1 & 8,1  & 9,1  &  \dots
&3,2 &   4,2 & 5,2  & 6,2  &  7,2 & 8,2  & 9,2  & \dots\\
\hline
  s_{2,1}   & -1  & 1  & 0 & 1  &  -1 & 0  &  \dots
 & 0 & -1   & 1  & 0  & 1  & -1  &0   & \dots\\
\hline
 s_{3,1}   & -1  & 0  & 1 & 1  &  0 & -1  &  \dots
 & -1 & -1   & 0  & 1  & 1  & 0  &-1   & \dots \\
\hline
\end{array}
$$
}

We have  periodicity after $6$ places. Further one has to cut after
the indices $(m+2,1)$ and  $(m+2,2)$.

\vspace{0,2cm}
 (a2)
Computation of  $\Hom(T_2,M)$:

 A homomorphism $f: T_2 \ra M$ is given by matrices
$P =(p_i) \in\matring_{m+2,1}(K)$,
$Q' =(q'_i) \in\matring_{m+1,1}(K)$,
$Q =(q_i) \in\matring_{m+2,1}(K)$,
 and
$S=(s_{i})\in \matring_{m+3,1}(K)$
such that
$S  =  X^{m+2}_{m+3} P$,
$Q  =  Y^{m+1}_{m+2} Q'$ and
$S  =  Y^{m+2}_{m+3} Q$.

It is easy to see that this yields the vanishing
conditions $s_1=0$,  $s_2=0$   and $s_{m+3}=0$.
We obtain  $\dim_K\Hom(T_2,M)=m $
and  a basis is given by $(m+3) \times 1$-matrices
$w^{(2)}_{3},w^{(2)}_{4}, \dots w^{(2)}_{m+2}$
where $w^{(2)}_{i}$ is the matrix with entries
$s_i=1$ and $s_j=0 $ for $j \neq i$.

\vspace{0,2cm}
 (a3)
Computation of  $\Hom(T_3,M)$:

A homomorphism $f: T_3 \ra M$ is given by matrices
$P' =(p'_i) \in\matring_{m+1,1}(K)$,
$P =(p_{i,j}) \in\matring_{m+2,2}(K)$,
$Q =(q_i) \in\matring_{m+2,1}(K)$,
$R=(r_i) \in\matring_{m+1,1}(K)$ and
$S=(s_{i})\in \matring_{m+3,2}(K)$
such that
$ P  X^{1}_{2}=  X^{m+1}_{m+2} P' $,
$S  =  X^{m+2}_{m+3}    P$,
$S     Y^{1}_{2}       =  Y^{m+2}_{m+3}    Q$,
 and
$S Z^{1}_{2}=Z^{m+1}_{m+3}R$.

The first equations yields $s_{1,2}=0$,
 $s_{m+2,1}=0$,  $s_{m+3,1}=0$ and   $s_{m+3,1}=0$.
Furthermore the last equation implies, similarly as in
(a2) that two further entries, say $s_{1,1}$ and
 $s_{2,1}$ are linearly dependent of the others coefficients,
however the precise formula for that will not be needed.
We obtain
  $\dim_K\Hom(T_3,M)=2m$ and a basis is given by
 $(m+3) \times 2$-matrices
$w^{(3)}_{3,1}, w^{(3)}_{4,1}   \dots,w^{(3)}_{m+1,1},
 w^{(3)}_{2,2},  w^{(3)}_{3,2}   \dots,w^{(3)}_{m+2,1}$
where $w^{(3)}_{i,j}$ is the matrix with entries
$s_{i,j }=1  $,  all other entries are $0$ except maybe
 $s_{1,1}$ and  $s_{2,1}$.

\vspace{0,2cm}
 (a4)
Computation of  $\Hom(T_4,M)$:

 A homomorphism $f: T_4 \ra M$ is given by matrices
$P' =(p'_i) \in\matring_{m+1,1}(K)$,
$P =(p_i) \in\matring_{m+2,1}(K)$,
$Q =(q_i) \in\matring_{m+2,1}(K)$,
 and
$S=(s_{i})\in \matring_{m+3,1}(K)$
such that
$P  =  X^{m+1}_{m+2} P'$,
$S  =  X^{m+2}_{m+3} P$,
 and
$S  =  Y^{m+2}_{m+3} Q$.

Again it is easy to verify that we get the following
vanishing conditions
 $s_{1}=0$,  $s_{m+2}=0$    and  $s_{m+3}=0$.
Thus we have  $\dim_K\Hom(T_4,M)=m $
and  a basis is given by $(m+3) \times 1$-matrices
$w^{(4}_{2},w^{(3)}_{4}, \dots w^{(4)}_{m+1}$
where $w^{(4)}_{i}$ is the matrix with entries
$s_i=1$ and $s_j=0 $ for $j \neq i$.

\vspace{0,2cm}
 (a5)
Computation of  $\Hom(T_5,M)$:

 A homomorphism $f: T_5 \ra M$ is given by matrices
$P' =(p'_i) \in\matring_{m+1,1}(K)$
$P =(p_i) \in\matring_{m+2,1}(K)$,
$Q'=(q'_i) \in\matring_{m+1,1}(K)$,
$Q =(q_i) \in\matring_{m+2,1}(K)$,
$R=(r_i) \in\matring_{m+1,1}(K)$ and
$S=(s_{i,j})\in \matring_{m+3,2}(K)$
such that
$ P  =  X^{m+1}_{m+2}    P' $,
$S  X^{1}_{2}=  X^{m+2}_{m+3}    P$,
$Q  =   Y^{m+1}_{m+2}  Q'$,
$S  Y^{1}_{2}  =  Y^{m+2}_{m+3}   Q$,
 and
$SZ^{1}_{2}=Z^{m+1}_{m+3}R$.

>From the first equations we infer that
 $s_{1,2}=0$,  $s_{2,2}=0$,
 $s_{m+2,1}=0$ and  $s_{m+3,1}=0$.
Moreover,  the last equation implies, similarly as in
(a2) that two further entries, say $s_{1,1}$ and
 $s_{2,1}$ are linearly dependent of the others.
We obtain
  $\dim_K\Hom(T_5,M)=2m$ and a basis is given by
 $(m+3) \times 2$-matrices
$w^{(3)}_{3,1}, w^{(3)}_{4,1}   \dots, \\
 w^{(3)}_{m+1,1},
 w^{(3)}_{3,2},  w^{(3)}_{3,2}   \dots,w^{(3)}_{m+3,2}$
where $w^{(3)}_{i,j}$ is the matrix with entries
$s_{i,j }=1  $,  all   other entries are $0$, except maybe
 $s_{1,1}$ and  $s_{2,1}$.

\vspace{0,2cm}
 (a6)
Computation of  $\Hom(T_6,M)$:

 A homomorphism $f: T_6 \ra M$ is given by matrices

$T=(t_i) \in\matring_{m,1}(K)$,
$P' =(p'_i) \in\matring_{m+1,1}(K)$
$P =(p_i) \in\matring_{m+2,1}(K)$,
$Q'=(q'_i) \in\matring_{m+1,1}(K)$,
$Q =(q_i) \in\matring_{m+2,1}(K)$,
$R=(r_i) \in\matring_{m+1,1}(K)$ and
$S=(s_{i,j})\in \matring_{m+3,2}(K)$
such that
$ P'  =  X^{m}_{m+1}    T $,
$ P  =  X^{m+1}_{m+2}    P' $,
$S  X^{1}_{2}=  X^{m+2}_{m+3}    P$,
$ Q'  =  Y^{m}_{m+1}    T $,
$Q  =   Y^{m+1}_{m+2}  Q'$
$S  Y^{1}_{2}  =  Y^{m+2}_{m+3}   Q$,
$R  =   Y^{m+1}_{m+2}  Q'$
$SZ^{1}_{2}=Z^{m+1}_{m+3}R$.

It is easily calculated that the equations imply that $S$ is of the
form

{\footnotesize
$$S=\begin{bmatrix}
t_1  &0 \\
t_2  &0 \\
t_3  &0 \\
t_4 &t_1 \\
\vdots & \vdots \\
t_m & t_{m-3} \\
0 & t_{m-2} \\
0 &  t_{m-1} \\
0 &  t_{m} \\
\end{bmatrix}
.$$
}

Therefore $\dim_K\Hom(T_{n+1},M)=m$ and a basis is given by $(m+2) \times
2$-matrices
$w^{(6)}_1  , w^{(6)}_2  , \dots,  w^{(6)}_m$  ,
where  $ w^{(6)}_i$ is the $(m+2)\times 2$-matrix  with  entries
$s_{i,1}=1$, $s_{i+3,2}=1$ and all other entries are zero.

We want to determine the matrices of the representation $N=F(M)$ of
the extended Dynkin quiver $\Gamma$ in the bases of the vector
spaces constructed above. For this reason we have to multiply with
the corresponding matrices $S^i_j$ from the right. In particular the
map $N(1) \ra N(0)$ is given by the formula

{\small
$$\begin{bmatrix}
0 &0 \\ s_{2,1} & 0  \\s_{3,1} &s_{3,2}  \\ \vdots \\ s_{m+2,1}&
s_{m+2,2}&  \\  0 &  s_{m+3,2}&\\
\end{bmatrix}
\mapsto
\begin{bmatrix}
0&0 &0 \\0& s_{2,1} & 0  \\0&s_{3,1} &s_{3,2}  \\ \vdots &\vdots \\0& s_{m+2,1}&
s_{m+2,2}&  \\ 0& 0 &  s_{m+3,2}&\\
\end{bmatrix}.
$$
}

This map in the bases
$w^{(1)}_{4,1}, w^{(1)}_{5,1}   \dots,w^{(1)}_{m+2,1},
 w^{(1)}_{3,2}, w^{(1)}_{4,2}   \dots,w^{(1)}_{m+3,1}$
and
$w^{(0)}_{3,1}, w^{(0)}_{4,1},  \\
 \dots,w^{(0)}_{m+1,1},
w^{(0)}_{2,2}, w^{(0)}_{3,2}, \dots,w^{(0)}_{m+2,2},
w^{(0)}_{3,3}, w^{(0)}_{4,3}, \dots,w^{(0)}_{m+3,3},
$
has the following form

{\footnotesize
$$ A=
\begin{array}{  |   ccc ccc    c ||    c|  ccc    ccc   c   |  }
\hline
& & &   & & &    & &     & & &  & & & \\
& & &   & & &    & &     & & &  & & & \\
& & &   & & &    & &     & & &  & & & \\
& & &   & & &    & &     & & &  & & & \\
& & &   & & &    & &     & & &  & & & \\
& & &   & & &    & &     & & &  & & & \\
& & &   & & &    & &     & & &  & & & \\
\hline
-1&1 &0 & 1  &-1 &0 & \dots   &0 &  -1  &1 & 0& 1 & -1&0 &\dots \\
-1&0 &1 & 1  &0 &-1 & \dots   &-1 &  -1   &0 &1 &1  &0 &-1 &\dots \\
1 & & &   & & &    & &     & & &  & & & \\
&1& &   & & &    & &     & & &  & & & \\
& &1&   & & &    & &     & & &  & & & \\
& & &1  & & &    & &     & & &  & & & \\
& & &   &1& &    & &     & & &  & & & \\
& & &   & &\ddots &    & &     & & &  & & & \\
\hline
& & &   & & &    &1 &     & & &  & & & \\
& & &   & & &    & & 1    & & &  & & & \\
& & &   & & &    & &     &1 & &  & & & \\
& & &   & & &    & &     & &1 &  & & & \\
& & &   & & &    & &     & & &1  & & & \\
& & &   & & &    & &     & & &  &\ddots & & \\
\hline
\end{array}
$$
} \noindent where the four indicated ``half'' rows with entries $0$,
$1$ and $-1$ have period $6$ but stop on the right hand side at the
position $m+2$ respectively  $m+3$.

Investigating the other linear maps in the same way we get the
following series of indecomposable representations $N_m$ for $Q$

   $$\xy\xymatrixcolsep{3.pc}\xymatrix{
  &     &        K^m      \ar @{->}[d]^-{ B  }   &      &  \\
         &      &  K^{2m}   \ar @{->}[d]^-{ A }    &    &     \\
    K^m    \ar @{->}[r]^-{ D  }  &    K^{2m} \ar @{->}[r]^-{ C  }    &       K^{3m+1}
   &     \ar @{->}[l]_-{E }    K^{2m}   &    \ar@{->}[l]_-{ G  }       K^m   }
  \endxy
  $$

where $A$ is defined above and

$$
B= \boxed{\begin{array}{c}
    \mathbf{0} \mid  \mathbf{ I_{m-1}}\\
    \hline
     \mathbf{- I_m}\\
    \hline
    0\dots 0\\
\end{array}}  \quad
C=
\boxed{\begin{array}{c|c}
  \mathbf{I_{m-1} }& \mathbf{0} \\
\hline
 \mathbf{0} &  \mathbf{I_{m+1} }\\
\hline
 \mathbf{0} & \mathbf{0} \\
\end{array}}
 \quad
D= \boxed{\begin{array}{c}
 \mathbf{0} \mid  \mathbf{I_{m-1}}\\
 \hline
  \mathbf{-I_m}\\
 \hline
 0\dots 0\\
\end{array}}
$$

$$
E=
\boxed{\begin{array}{c|c}
 \mathbf{ I_{m-1}} & \mathbf{0} \\
\hline
 \mathbf{0} & \mathbf{0}  \\
\hline
 \mathbf{0} &     \mathbf{I_{m+1}} \\
\end{array}}
\quad
G=
\boxed{\begin{array}{c}
    \mathbf{0} \mid  \mathbf{I_{m-2}}\\
 \hline
 0\dots 0\\
\hline
 0\dots 0\\
 \hline
  \mathbf{I_m}
\end{array}}
$$

Applying the functor $F$  to  the second series of rank $3$ modules
over $\La$ the same method yields the following indecomposable
representations for $\Gamma$.

   $$\xy\xymatrixcolsep{3.pc}\xymatrix{
  &     &        K^m      \ar @{->}[d]^-{ B  }   &      &  \\
         &      &  K^{2m+1}   \ar @{->}[d]^-{ A }    &    &     \\
    K^m    \ar @{->}[r]^-{ D  }  &    K^{2m+1} \ar @{->}[r]^-{ C  }    &       K^{3m+1}
   &     \ar @{->}[l]_-{E }    K^{2m+1}   &    \ar@{->}[l]_-{ G  }
   K^m   }
  \endxy
  $$

where

$$
A= \boxed{\begin{array}{c|c}
     \mathbf{I_m }& \mathbf{0}\\
    \hline
    \mathbf{0}& \mathbf{ I_{m+1}}\\
    \hline
    \mathbf{0} & \mathbf{0}\\
  \end{array}}  \quad
   B= \boxed{\begin{array}{c}
    \mathbf{ I_m}\\
    \hline
     \mathbf{-I_m}\\
    \hline
    0\dots 0\\
\end{array}}  \quad
C=
\boxed{\begin{array}{c|c}
  \mathbf{I_{m}} & \mathbf{0} \\
\hline
 \mathbf{0} & \mathbf{0} \\
\hline
 \mathbf{0} &  \mathbf{I_{m+1}} \\
\end{array}}$$

$$D= \boxed{\begin{array}{c}
 \mathbf{0} \mid  \mathbf{I_{m-1}}\\
 \hline
  0\dots 0\\
  0\dots 0\\
  \hline
   \mathbf{ I_m}
\end{array}}\quad E=
\boxed{\begin{array}{c|c}
    \mathbf{0} & \mathbf{0} \\
    \hline
    1\ -1\ 1\ -1\dots & 1\ -1\ 1\ -1\dots\\
    \hline
     \mathbf{I_m }& \mathbf{0}  \\
    \hline
    \mathbf{0} &    \mathbf{ I_{m+1} }
\end{array}}
\quad
G=
\boxed{\begin{array}{c}
    \mathbf{I_m}\\
 \hline
    \mathbf{-I_m}\\
\hline
 0\dots 0
\end{array}}
$$

{\footnotesize
\bibliographystyle{amsplain}

\providecommand{\bysame}{\leavevmode\hbox to3em{\hrulefill}\thinspace}
\providecommand{\MR}{\relax\ifhmode\unskip\space\fi MR }
\providecommand{\MRhref}[2]{%
  \href{http://www.ams.org/mathscinet-getitem?mr=#1}{#2}
}
\providecommand{\href}[2]{#2}

\end{document}